\documentclass{amsart}
\usepackage{graphicx}
\newtheorem{thm}{Theorem}

\theoremstyle{definition}

\theoremstyle{remark}

\newcommand{\nbd}{neighborhood}

\begin{document}
\title{Non-vanishing complex vector fields and the Euler
  characteristic}
\thanks{Submitted for publication July 25, 2008.}
\author{Howard Jacobowitz}
\address{Rutgers University\\
Camden, New Jersey 08012}
\email{jacobowi@camden.rutgers.edu}
\begin{abstract}
Every manifold admits a nowhere vanishing complex vector field.  If,
however, the manifold is compact and orientable and the complex
bilinear form associated to a Riemannian metric is never zero when
evaluated on the vector field, then the manifold must have zero Euler 
characteristic.
\end{abstract}
\maketitle
One of the oldest and most basic results in global differential
topology relates the topology of a manifold to the zeros of its
vector fields. Let $M$ be a compact and orientable manifold and let
$\chi (M)$ denote its Euler characteristic.  Here is the simplest
statement of this relation.
\begin{eqnarray}\label{real}
\mbox{If there is a global nowhere zero vector field on}\  M \mbox{
then}\  \chi (M) =0.
\end{eqnarray}
This of course is for a real vector field  (that is, for a section
$M\to TM$).  On the other hand, it is easy to see that any manifold
admits a nowhere zero complex vector field. (A complex vector field is a section
$M\to \mathbb{C}\otimes TM$). This can be seen most simply by observing
that a generic perturbation of any section, even the zero section
itself, must be everywhere different from zero.

It is natural to seek a condition on a nowhere zero
complex vector field which would again imply $\chi (M) =0$. Curiously,
a trivial restatement of (\ref{real})  leads to such a condition. Let $g$ be
any Riemannian metric on $M$.

\bigskip
\noindent (2)\
 Let  $v : M\to TM $ be
a global vector field on $M$. If the Riemannian metric $g(v,v)$ is
never  zero, then $\chi (M)=0$.
\bigskip

\noindent Here is the condition for complex vector fields.
\begin{thm}\nonumber
Let $v : M\to \mathbb{C}\otimes TM$  be a global vector field on $
M$. If the bilinear form  g(v,v) is never zero, then $ \chi (M)=0$.
\end{thm}

\noindent Here $g$ is extended to complex vector fields by taking
$g(v,w)$ to be complex linear in each argument;  for $v=\xi + i\eta$ we have
\[
g(v,v)=g(\xi ,\xi )-g(\eta ,\eta )+2ig(\xi ,\eta ).
\]
\begin{proof}
We show that if $g(v,v)\not= 0$ then $v$ can be deformed to a
nowhere zero real vector field. So the Euler characteristic would be
zero, according to  (\ref{real}).   We decompose $M$ as
\[
M=A_+\cup B\cup A_-
\]
where $g(\xi ,\xi )>g(\eta ,\eta )$ on $A_+$, the opposite
inequality holds on $A_-$, and equality holds on $B$.  We assume for
now that $B$ is not empty.  Note that $\xi$ is nowhere zero in $A_+$
and $\eta$ is nowhere zero in $A_-$. Further, since $g(v,v)$ is
never zero, we have that $g(\xi ,\eta )$ is never zero on $B$.  Thus
there is an open \nbd\ $\Omega$ of $B$ on which $g(\xi ,\eta )$ is
never zero. We may take $\Omega$ to have a smooth boundary.  We have
that $\eta$ is never zero in $A_-\cup\Omega$ and $\xi$ is never zero
in $A_+\cup\Omega$. Let $\Omega _1$ be an open set chosen so that
\[
B\subset \Omega _1,\mbox{      } \overline{\Omega _1}\subset \Omega
\]
and
\item
\[\overline{\Omega _1} \mbox{  is a \nbd\  retract of }\overline{
\Omega}.
\]

The boundary of $\Omega$ has two components, one in $A_+$ and the
other in $A_-$. (That is, the boundary of $\Omega$ is the union of two
sets, neither of which need be connected.)  The same is true for the boundary of $\Omega _1$.
We will work only with the components in $A_+$.  Call them $\Sigma$
and $\Sigma _1$.  Each of these sets separates $M$ into two
components.  We seek to deform $v$ to a nowhere vanishing real
vector field $u$.  Set $u=\xi$ on the component of $M-\Sigma$ which
does not contain $A_-$.  The sets $\Sigma$ and $\Sigma _1$ bound a
region which retracts onto $\Sigma _1$.  We want to rotate $\xi$ to
$\eta$, (or to $-\eta$) as the retraction takes $\Sigma$ to $\Sigma
_1$. Since $g(\xi ,\eta ) \neq 0$, in this region, this is easily
done. Pick a point in this region.  The angle $\theta$ between the
vectors $\xi$ and $\eta$ satisfies one of the alternatives
\[
0\leq \theta \leq \pi /2 \mbox{  or  } \pi /2 < \theta \leq \pi
\]
and whichever alternative is satisfied at that point is also
satisfied at all points in the region.  Thus as we retract $\Sigma$
to $\Sigma _1$, we may rotate $\xi$ to $\eta$, or, respectively to 
$-\eta$. Finally, define $u=\eta$, respectively $u= -\eta$, on the
component of $M-\Sigma _1$ which contains $A_-$.

If $B$ is empty, the proof is even
easier.  Now either $g(\xi ,\xi)>g(\eta ,\eta )$ everywhere and so
$\xi$ is a nowhere zero real vector field or the opposite
inequality holds and $\eta$ is a nowhere zero real vector field.
\end{proof}
{\it{Remark}}.  \noindent We have proved the theorem by reducing to
(\ref{real}).  This latter result goes back to H. Hopf; an
influential modern proof was given by Atiyah \cite{A}.  Atiyah's
proof makes use of the Clifford algebra structure on the bundle of
exterior forms. Our Theorem can be proved directly, without reducing
to (\ref{real}), by following Atiyah's proof using the corresponding
complex Clifford algebra.

There is a stronger version of (\ref{real}) which expresses the
Euler characteristic as the algebraic sum of the indices of the
zeros of the vector field.  (Indeed this is the result of Hopf.)  It
would be interesting to generalize this to complex vector fields.

\end{document}